\documentclass{amsart}
\usepackage{amsmath}
\usepackage{amsthm}
\usepackage{amssymb}
\usepackage{epsfig}

\newcommand{\n}{{\mathcal{N}_2}}
\newcommand{\C}{{\mathcal{C}}}
\newcommand{\bm}[1]{\mbox{\boldmath $#1$}}

\theoremstyle{plain}
\newtheorem{thm}{Theorem}
\newtheorem{lem}[thm]{Lemma}
\newtheorem*{tc}{Technical Condition}
\theoremstyle{definition}
\newtheorem*{example}{Example}
\newtheorem*{defn}{Definition}

\title{A Non-Messing-Up Phenomenon for Posets}
\author{Bridget Eileen Tenner}
\date{October 17, 2005}

\begin{document}

\begin{abstract} We classify finite posets with a particular sorting
property, generalizing a result for rectangular arrays. Each poset is
covered by two sets of disjoint saturated chains such that, for any
original labeling, after sorting the labels along both sets of chains, the
labels of the chains in the first set remain sorted. We also characterize
posets with more restrictive sorting properties. \end{abstract}

\maketitle

\section{Introduction}

The so-called Non-Messing-Up Theorem is a well known sorting result for
rectangular arrays.  In \cite{knuth}, Donald E.~Knuth
attributes the result to Hermann Boerner, who mentions it in a footnote in
Chapter V, \S 5 of \cite{boerner}.  Later, David Gale and Richard M.~Karp
include the phenomenon in \cite{galepreprint} and in
\cite{gale}, where they prove more general results about order
preservation in sorting procedures.  The first use of the term
``non-messing-up'' seems to be due to Gale and Karp, as suggested in
\cite{knuthletter}.  One statement of the result is as follows.

\begin{thm}
\label{nmu}
Let $A = (a_{ij})$ be an $m$-by-$n$ array of real numbers.  
Put each row of $A$ into non-decreasing order.  That is, for each
$1 \le i \le m$, place the values $\{a_{i1}, \ldots, a_{in}\}$ in
non-decreasing order (henceforth denoted \emph{row-sort}).  This yields the array $A' = (a_{ij}')$.  Column-sort $A'$.  Each row in the resulting array is in non-decreasing order.
\end{thm}

\begin{proof}
Following the solution submitted by J.~L.~Pietenpol in \cite{monthly}, first row-sort $A$ to obtain $A'$.  Column-sort $A'$ by first permuting the rows to order the first column; to order the second column, permute the rows without their first entries; and so on.  The rows remain sorted at each step of the procedure.
\end{proof}

Applying the theorem to the transpose of the array $A$, the
sorting can also be done first in the columns, then in the rows, and the
columns remain sorted.  

\begin{example}
\[\begin{matrix}4 & 9 & 7 & 8\\12 & 5 & 1 & 10\\2 & 6 & 11 & 3\end{matrix}
\ \ \xrightarrow{\text{row-sort}} \ \ \begin{matrix}4 & 7 & 8 & 9\\1 & 5 & 10 &
12\\2 & 3 & 6 &11\end{matrix} \ \ \xrightarrow{\text{column-sort}} \ \
\begin{matrix}1 & 3 & 6 & 9\\2 & 5 & 8 & 11\\4 & 7 & 10 & 12\end{matrix}\]
\end{example}

Answering a question posed by Richard P.~Stanley, the author's thesis
advisor, this paper defines a notion of non-messing-up for posets and Theorem~\ref{answer} generalizes Theorem~\ref{nmu} by characterizing all posets with this property.

\medskip

Standard terminology from the theory of partially ordered sets will be used throughout the paper.  A good reference for these terms and other information about posets is Chapter 3 of \cite{ec1}.

The rectangular array in Theorem~\ref{nmu} can be viewed as the poset $\bm{m} \times \bm{n}$ (where $\bm{j}$ denotes a $j$-element chain).  The \emph{rows} and \emph{columns} are two different sets of disjoint saturated chains, each covering this poset.  Sorting a chain orders the chain's labels so that the minimum element in the chain has the minimum label.  Thus, sorting the labels $\{1, \ldots, mn\}$ in this manner gives a linear extension of $\bm{m} \times \bm{n}$.

\begin{figure}[htbp]
\centering
\epsfig{file=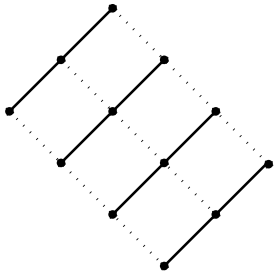}
\caption{The poset $\bm{3} \times \bm{4}$ with chain covers marked in solid and dotted lines.}
\label{4x3grid}
\end{figure}

\begin{defn}
An \emph{edge} in a poset $P$ is a covering relation $x \lessdot y$.  Two elements in $P$ are \emph{adjacent} if there is an edge between them.
\end{defn}

\begin{defn}
A \emph{chain cover} of a poset $P$ is a set of disjoint saturated chains covering the elements of $P$.
\end{defn}

\begin{defn}
A finite poset $P$ has the \emph{non-messing-up} property if there exists an unordered pair of chain covers $\{\C_1, \C_2\}$ such that
\begin{enumerate}
\item For any labeling of the elements of $P$, $\C_i$-sorting and then $\C_{3-i}$-sorting leaves the labels sorted along the chains of $\C_i$, for $i = 1$ and $2$; and
\item Every edge in $P$ is contained in an element of $\C_1$ or $\C_2$.
\end{enumerate}
\noindent The set $\n$ consists of all posets with the non-messing-up property, where the subscript indicates that an unordered \emph{pair} of chain covers is required.  For a non-messing-up poset $P$ with chain covers as defined, write $P \in \n$ \emph{via} $\{\C_1, \C_2\}$.
\end{defn}

The map $\lambda: P \rightarrow \mathbb{R}$ will denote the labels of elements of a poset $P$.  The main result of this paper is the classification of $\n$.

It is illustrative to clarify the difference between this result and Gale and Karp's work in \cite{galepreprint} and \cite{gale}.  Gale and Karp consider a poset $P$ and a partition $F$ of the elements of $P$.  The elements in each block of $F$ are linearly ordered, not necessarily in relation to comparability in $P$.  Given $P$ and $F$, the authors determine whether each natural labeling of $P$, sorted within each block of $F$, yields a labeling that is still natural.  This paper does not require that the original labeling be natural.  In fact, it is the labelings that are not natural and that do not become natural after the first sort that determine membership in $\n$.  Additionally, the partition blocks in $\n$ are saturated chains, and every covering relation must be in at least one of these chains.  The goal of this paper is to determine, for a given poset, when there \emph{exist} chain covers with the non-messing-up property, not if a given pair of chain covers has the property.

It is important to emphasize that $\{\C_1,\C_2\}$ is an unordered pair and that there is a symmetry between the chain covers.  The elements of $\C_1$ and their edges will be referred to as \emph{red}, and the elements of $\C_2$ and their edges as \emph{blue}.   If an edge belongs to both chain covers, it is \emph{doubly colored}.  The symmetry between the chain covers may be expressed by a statement about red and blue chains and an indication that a \emph{color reversed} version of the statement is also true.

Central to the classification of $\n$ is the following object, defined and denoted $\mathcal{C}_{kn}$ in \cite{postnikov}.  The notation is changed here to avoid confusion with the chain covers.

\begin{defn}
Fix positive integers $k$ and $n$ so that $k < n$.  The \emph{cylinder poset} $Cyl_{k,n}$ is $\mathbb{Z}^2$ modulo the equivalence relation $(i,j) \sim (i-k, j+n-k)$.  That is,
\[Cyl_{k,n} = \mathbb{Z}^2/(-k,n-k)\mathbb{Z}.\]
\noindent The partial order on $Cyl_{k,n}$ is induced by the componentwise partial order on $\mathbb{Z}^2$.
\end{defn}

Note that the product of two chains is always a convex subposet of a cylinder poset.

The classification in Theorem~\ref{answer} states that $\n$ is the set of disjoint unions of finite connected posets that each can be ``reduced'' to
a convex subposet of a cylinder poset,
subject to a technical constraint.  Informally speaking, $P$ reduces to
$Q$ if $P$ is formed by replacing particular elements of $Q$ with chains
of various lengths.  Sample Hasse diagrams for elements of $\n$ are shown in Figures~\ref{treeanswer}(a), \ref{bigcircuit}(a), \ref{typeiistrings}(a),
and~\ref{bigcrown}(a).

Section~2 of this paper addresses definitions and preliminary results.  The definitions describe the objects and operations needed for the classification, and the results will be the fundamental tools for defining $\n$.  The main theorem is proved in Section~3 by induction on the size of a connected poset.  The final section of the paper discusses further directions for the study of non-messing-up posets, including several open questions.

\section{Preliminary results}

The definition of a non-messing-up poset requires that every edge be colored.  Therefore, as in the case of the product of two chains, $\C_i$-sorting any labeling of a poset $P$ and then $\C_{3-i}$-sorting yields a linear extension of $P$ if the labels are $\{1, \ldots, |P|\}$.  The chains of $\C_i$ are disjoint, so each element of a non-messing-up poset is covered by at most two elements, and covers at most two elements.

It is sufficient to consider connected posets, as a poset is in $\n$ if and only if each of its connected components is in $\n$.  Key to determining membership in $\n$ is the following fact.

\begin{thm}\label{convexity}
Every convex subposet of an element of $\n$ is also in $\n$.
\end{thm}

\begin{proof}
Consider $P \in \n$ via $\{\C_1,\C_2\}$.  Let $Q$ be a convex subposet of $P$.  Consider a labeling $\lambda$ of $Q$, and let $m = \min_{x \in Q} \lambda(x)$ and $M = \max_{x \in Q} \lambda(x)$.  Extend $\lambda$ to a labeling $\widetilde{\lambda}$ of $P$ by
\[\widetilde{\lambda}(x) = \left\{ \begin{array} {c@{\quad:\quad}l} \lambda(x) & x \in Q; \\ m & x \notin Q \text{ and } x < y \text{ for some } y \in Q; \\ M & x \notin Q \text{ and } x > y \text{ for some } y \in Q; \\ m & \text{otherwise.} \end{array} \right.\]
\noindent The convexity of $Q$ makes this well-defined.  For the labeling $\widetilde{\lambda}$ of $P$, $\C_i$-sort and then $\C_{3-i}$-sort.  As $P \in \n$, the labels remain sorted along the chains of $\C_i$.  By construction, the only elements whose labels may change during sorting are in the subposet $Q$.  Thus $Q \in \n$ via $\{\C_1\vert_Q, \C_2\vert_Q\}$.
\end{proof}

Consider a poset with chain covers $\{\C_1,\C_2\}$.  For any labeling, $\C_i$-sort and then $\C_{3-i}$-sort.
The $\C_i$ chains are still in order if the labels on each edge of each chain in $\C_i$ are in order.  Some of these edges will automatically have sorted labels: for example, doubly colored edges.  Additionally, suppose that $x \lessdot y$ is an edge in a chain $\bm{c} \in \C_i$, and $x$ and $y$ are both in $\C_{3-i}$ chains that are entirely contained within $\bm{c}$.  In this situation, $\lambda(x)$ and $\lambda(y)$ will be unchanged after the $\C_{3-i}$-sort, and so will necessarily remain sorted.

\begin{lem}\label{chains}
If a convex subposet of a non-messing-up poset is a chain, then there is a 
red chain or a blue chain containing this entire subposet. \end{lem}

\begin{proof}
Consider $P \in \n$ via $\{\C_1,\C_2\}$.  Let $\bm{c}$ be a convex chain of $P$.  Suppose neither $\C_1$ nor $\C_2$ has an element containing $\bm{c}$.  Then there are chains $\bm{c}_i \in
\C_i$ where, up to color reversal, $\widetilde{\bm{c}}_1$ extends
below $\widetilde{\bm{c}}_2$, and $\widetilde{\bm{c}}_2$
extends above $\widetilde{\bm{c}}_1$ (for $\widetilde{\bm{c}}_i = \bm{c}_i \cap \bm{c}$).  Consider $\widetilde{\bm{c}}_1
\cup \widetilde{\bm{c}}_2 = \{x_1 \lessdot \cdots \lessdot x_k\}$.  Set 
\[\lambda(x_1, \ldots, x_k) = (2, \ldots, |\widetilde{\bm{c}}_1| + 1, 1, |\widetilde{\bm{c}}_1| + 2, \ldots, k).\]
\noindent The $\C_1$-sort will not change any labels, but the $\C_2$-sort will label a non-minimal element of $\widetilde{\bm{c}}_1$ with $1$, a contradiction.
\end{proof}

\begin{defn}
A \emph{diamond} in a poset is a convex subposet that is the union of distinct (saturated) chains which only intersect at a common minimal element and a common maximal element.
\end{defn}

\begin{lem}\label{diamondcolors}
Let $Q$ be a diamond consisting of chains $\bm{a}$ and $\bm{b}$ in a non-messing-up poset.  Let $x$ be the minimal element in $\bm{a}$ and $\bm{b}$, denoted $\min(\bm{a})$ and $\min(\bm{b})$, and let $y = \max(\bm{a}) = \max(\bm{b})$, with similar notation.  Up to color reversal, one of the following is true (where $\bm{c}\setminus z$ is taken to mean $\bm{c}\setminus\{z\}$).
\begin{enumerate}
\item There exists a red chain containing $\bm{a}\setminus y$, a blue chain containing $\bm{b}\setminus y$, a red chain containing $\bm{b}\setminus x$ and a blue chain containing $\bm{a}\setminus x$; or
\item There exists a red chain containing $\bm{a}$ and a blue chain containing $\bm{b}$.
\end{enumerate}
\noindent Call the former of these a \emph{Type I} diamond and the latter a \emph{Type II} diamond.
\end{lem}

\begin{proof}
Let $\bm{a} = \{x = a_0 \lessdot \cdots \lessdot a_M = y\}$ and $\bm{b} = \{x = b_0 \lessdot \cdots \lessdot b_N = y\}$.  Suppose $x \lessdot a_1$ is red and $x \lessdot b_1$ is blue.  If $a_{M-1} \lessdot y$ is blue and $b_{N-1} \lessdot y$ is red, then Lemma~\ref{chains} makes $Q$ a Type I diamond.  If, instead, $a_{M-1} \lessdot y$ is red and $b_{N-1} \lessdot y$ is blue, then Lemma~\ref{chains} requires $Q$ to have Type II.
\end{proof}

\begin{figure}[htbp]
\begin{center}
$\begin{array}{c@{\hspace{1in}}c}
\multicolumn{1}{l}{\mbox{\bf{(a)}}} & \multicolumn{1}{l}{\mbox{\bf{(b)}}}\\
[-.45cm]
\epsfig{file=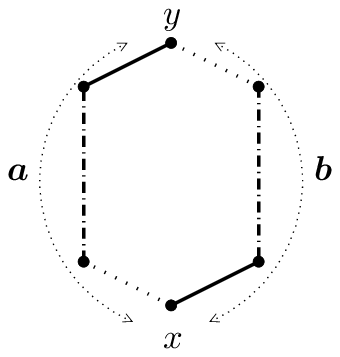}
& \epsfig{file=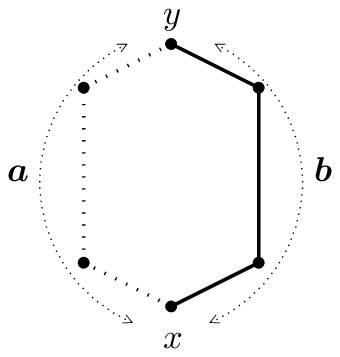}
\end{array}$
\end{center}
\caption{$\bf{(a)}$ Type I diamond coloring.  $\bf{(b)}$ Type II diamond coloring.  The intervals $\bm{a} \setminus \{x,y\}$ and $\bm{b} \setminus \{x,y\}$ may be partially or totally doubly colored.}
\label{diamondcoloring}
\end{figure}

\begin{defn} A diamond with \emph{bottom chain} of length $k$ and
\emph{top chain} of length $l$ is a convex subposet that is a diamond
with minimum $x$ and maximum $y$, a chain of $k$ elements covered by $x$,
and a chain of $l$ elements covering $y$, with no other elements or
relations among the elements already mentioned. \end{defn}

\begin{figure}[htbp]
\centering
\epsfig{file=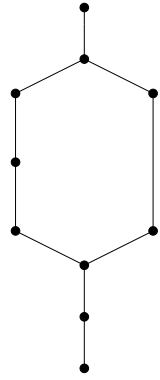}
\caption{A diamond with bottom chain of length $2$ and top chain of length $1$}
\end{figure}

\begin{defn}
A diamond in the cylinder poset $Cyl_{k,n}$ that is not a diamond in $\mathbb{Z}^2$ for any choice of preimages is said to \emph{go around} the cylinder.
\end{defn}

\begin{example}
In $Cyl_{3,4}$, the four element diamond consisting of the equivalence classes including $(0,1)$, $(1,1)$, $(1,0)$, and $(2,0)$ goes around the cylinder.
\end{example}

The technical condition mentioned in the introduction is due to the following requirement.

\begin{lem}\label{typeirules}
Let  $Q \subseteq P \in \n$ be a diamond consisting of chains $\bm{a}$ and $\bm{b}$.  Suppose there is a coloring of $P$ for which $Q$ has Type I, with bottom and top chains $\bm{C}$ and $\bm{D}$.  Then there are chains in that coloring such that, up to color reversal, $(\bm{C} \cup
\bm{a})\setminus y$ is red, $(\bm{C} \cup \bm{b})\setminus y$ is blue,
$(\bm{a} \cup \bm{D})\setminus x$ is blue, and $(\bm{b} \cup
\bm{D})\setminus x$ is red.  Also, $\max\{|\bm{C}|,|\bm{D}|\} <
\min\{|\bm{a}| - 2,|\bm{b}| - 2\}$. 
\end{lem}

\begin{proof}
Let $\bm{a} = \{x = a_0 \lessdot \cdots \lessdot a_M = y\}$, $\bm{b} = \{x = b_0 \lessdot \cdots \lessdot b_N = y\}$, $\bm{C} = \{w_k \lessdot \cdots \lessdot w_1\}$, and $\bm{D} = \{z_1 \lessdot \cdots \lessdot z_l\}$.  Color $Q$ so that it has Type I and $\bm{a} \setminus y$ is red.  By Lemma~\ref{chains}, $[w_k, a_{M-1}]$ is red, $[w_k, b_{N-1}]$ is blue, $[a_1, z_l]$ is blue, and $[b_1, z_l]$ is red.  Define $\lambda: Q \cup \bm{C} \cup \bm{D} \rightarrow \mathbb{R}$ by
\[\lambda(b_1, \ldots, b_N, z_1, \ldots, z_l, w_k, \ldots, w_1, a_0, \ldots, a_{M-1}) = (1, \ldots, M+N+k+l).\]
 
First $\C_1$-sorting and then $\C_2$-sorting, it is only necessary to check the labels of the edges $b_{N-1} \lessdot y$ and $x \lessdot a_1$.  The $\C_1$-sort changes nothing.  After the $\C_2$-sort,
\[\lambda(a_1) = N; \hspace{.075in} \lambda(b_{N-1}) = N+k+l+1; \hspace{.075in}
\lambda(x) = \left\{ \begin{array}{c@{\quad:\quad}l} k+1& k+2 \le
N; \\ k+l+2 & k+2 > N; \end{array}\right.\]
\[\text{and} \hspace{.075in}\lambda(y) = \left\{ \begin{array}{c@{\quad:\quad}l} M+N+k & l+2 \le M; \\ M+N-1 & l+2 > M. \end{array}\right.\]
\noindent The non-messing-up property requires $k+2 \le N$ and $l+2 \le M$.  An analogous color reversed argument yields $k+2 \le M$ and $l+2 \le N$.
\end{proof}

If $\max\{|\bm{C}|,|\bm{D}|\}=\min\{|\bm{a}| - 3, |\bm{b}| - 3\}$, then the
described chains have the non-messing-up property, so the bounds in Lemma~\ref{typeirules} are sharp.

As suggested by the main result, cylinder posets are crucial in the study of $\n$.

\begin{thm}\label{cylinder}
For all $k$ and $n$, any finite convex subposet of $Cyl_{k,n}$ is in $\n$.  The chain covers for this poset are of the same form as the chain covers in Theorem~\ref{nmu}.
\end{thm}

\begin{proof}
Let $P$ be a finite convex subposet of $Cyl_{k,n}$.  Cut $Cyl_{k,n}$ to get a preimage of $P$ in the plane.  Glue together copies of this poset via the identifications on the cylinder.  After perhaps removing elements at the edges of the planar poset, this is a convex subposet of $\bm{M} \times \bm{M}$ for some $M$.  For the labeling $\lambda: P \rightarrow \mathbb{R}$, label every preimage of $x$ in the plane by $\lambda(x)$.  Glue enough copies of the poset so that after the two sorts, the centermost copy in the plane has the labels it would have had on the cylinder.  This is possible because only finitely many elements cross over a line of identification.  Since $\bm{M} \times \bm{M} \in \n$, the labels of all the chains in the centermost copy of the cut poset must be in order.
\end{proof}

Before proving the main theorem, it remains to define reduction.  

\begin{defn}
The process of \emph{splitting} the element $x \in Q'$ gives a poset $Q$ where
\begin{enumerate}
\item $x \in Q'$ is replaced by $\{x_1 \lessdot
\cdots \lessdot x_{s(x)}\}$ for some positive integer $s(x)$;
\item All elements and relations in $Q'\setminus x$ are unchanged in $Q$;
\item If $y \gtrdot x$ in $Q'$, then $y \gtrdot x_{s(x)}$ in $Q$; and
\item If $y \lessdot x$ in $Q'$, then $y \lessdot x_1$ in $Q$.
\end{enumerate}
\noindent If $Q$ is formed by splitting elements of $\widetilde{Q}$, then $Q$ \emph{reduces} to $\widetilde{Q}$, denoted $Q \leadsto \widetilde{Q}$.
\end{defn}

\begin{defn}
Let $P \leadsto \widetilde{P} \in \n$.  The coloring of $\widetilde{P}$ \emph{induces} the coloring of $P$ if the edge $\widetilde{u} \lessdot \widetilde{v}$ in $\widetilde{P}$ and
its image, the edge $u \lessdot v$ in $P$, are colored in the same way.
Edges in the chain into which an element splits get doubly colored.
\end{defn}

\begin{figure}[htbp]
\begin{center}
$\begin{array}{c@{\hspace{.25in}}c@{\hspace{.25in}}c@{\hspace{.25in}}c}
\multicolumn{1}{l}{\mbox{\bf{(a)}}} & \multicolumn{1}{l}{\mbox{\bf{(b)}}} 
& \multicolumn{1}{l}{\mbox{\bf{(c)}}} & \multicolumn{1}{l}{\mbox{\bf{(d)}}} 
\\
[-.53cm]
\hspace{.35in}\epsfig{file=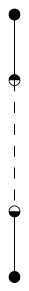} &
\hspace{.35in}\raisebox{.25in}{\epsfig{file=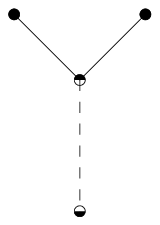}} &
\hspace{.35in}\epsfig{file=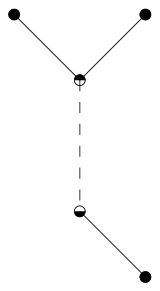} 
& \hspace{.35in}\epsfig{file=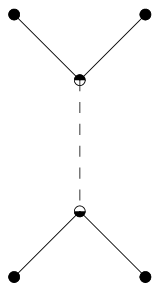}
\end{array}$
\end{center}
\caption{How to split different kinds of elements.}
\label{splits}
\end{figure}

\section{Characterization of $\n$}

The classification of the set $\n$ is done in two steps.  The first direction will show that any finite poset reducing to a convex subposet of a cylinder poset, subject to a technical constraint imposed by Lemma~\ref{typeirules}, has the non-messing-up property.  The second step will show the reverse inclusion.  Both directions are proved by induction on the size of a connected poset.

\begin{thm}\label{answer}
The collection $\n$ is exactly the set of finite posets each of whose connected components $P$ reduces to $\widetilde{P}$, a convex subposet of a cylinder poset, given the following stipulation:
\begin{center}
\parbox{4.125in}{\begin{tc} For any diamond $\{w \lessdot x,y \lessdot z\}$ in $\widetilde{P}$ that does not go around the cylinder,
\begin{center}
$\max\{s(w),s(z)\} \le \min\{s(x),s(y)\}$.
\end{center}
\end{tc}}
\end{center}
\noindent The required coloring of the connected poset $P \in \n$ is induced by the coloring of  $\widetilde{P}$, which is inherited from the coloring in Theorem~\ref{cylinder}.  
\end{thm}

\begin{proof}
The theorem is proved by induction on the size of a connected poset.  All one, two, and three element connected posets clearly have the non-messing-up property with coloring inherited from the coloring defined Theorem~\ref{nmu}, and they are all convex subposets of $\bm{3} \times \bm{3}$. These posets do not contain any diamonds, so the technical condition is trivially satisfied.  Assume inductively that the theorem holds for all connected posets with less than $K$ elements for $K \ge 4$.  Throughout the proof, let $P$ be a $K$ element connected poset.

Both directions of the proof consider a subposet $P'$ formed by removing 
either a maximal or a minimal element from $P$.  Thus $P'$ is convex in $P$, 
and it is not hard to see that the suppositions for $P$ must hold for $P'$ as 
well.  Each connected component in $P'$ has fewer than $K$ elements, so the 
theorem holds for $P'$ by the inductive assumption.

First let $P$ be a connected poset such that $P \leadsto \widetilde{P}$, a convex subposet of a cylinder poset, subject to the technical condition in the statement of the theorem.

Suppose that there is a maximal or minimal element $z \in P$ adjacent to only one element $v \in P$, and that $z$ and $v$ both map to $\widetilde{v}$ in $\widetilde{P}$.  Let $P' = P \setminus z$.  The technical condition holds for the poset $P'$ since removing a maximal or minimal element can at worst shorten the length of the top or bottom chain of a diamond.  Also, $P'\leadsto \widetilde{P}$.  The poset $P'$ has $K-1$ elements, so $P' \in \n$ via $\{\C_1,\C_2\}$ by the inductive hypothesis, and this coloring must be induced by the coloring of $\widetilde{P}$.  If $v \in \bm{c}_i \in \C_i$, let $\bm{c}_i^+ = \bm{c}_i \cup z$ and $\C_i^+ = (\C_i \setminus \bm{c}_i) \cup \bm{c}_i^+$.  Then $P \in \n$ via $\{\C_1^+,\C_2^+\}$.  Moreover, the coloring of $P$ must be induced by the coloring of its reduced poset, otherwise there would be a coloring of $P'$ contradicting this aspect of the inductive hypothesis.

If, on the other hand, there is no such $z \in P$, then either $P$ is a chain or there is a maximal or minimal element $w \in P$ adjacent to distinct elements.  A chain, itself a convex subposet of a cylinder poset, is in $\n$ via itself and any chain cover, as required by Lemma~\ref{chains}.  If $P$ is not a chain, let $w$, without loss of generality, be minimal and covered by $x$ and $y$.  Let $P' = P \setminus w$.  This poset reduces to a convex subposet $\widetilde{P}'$ of $\widetilde{P}$, so $P' \in \n$ via $\{\C_1,\C_2\}$ by the inductive hypothesis.

The coloring of $P'$ is induced by that of its reduced poset $\widetilde{P}' 
\subseteq \widetilde{P}$.  Thus there are $\bm{c}_i \in \C_i$ such that, 
without loss of generality, $\min(\bm{c}_1) = x$ and $\min(\bm{c}_2) =y$.  Let $\bm{c}_i^+ = \bm{c}_i \cup w$ and $\C_i^+ = (\C_i \setminus 
\bm{c}_i) \cup \bm{c}_i^+$.  Since $P' \in \n$, membership of $P$ in $\n$ 
depends on the labels of the edges $w \lessdot x$ and $w \lessdot y$ after the 
two sorts.  After $\C_1^+$-sorting, $\lambda(w) 
\le \lambda(x)$.  Subsequently $\C_2^+$-sorting cannot increase $\lambda(w)$, 
and $\lambda(w) \le \lambda(y)$.  Suppose $x \in \bm{c}_2' \in \C_2^+$.  If 
the label of $x$ decreases after the second sort, then the set $S = \{ v \in 
\bm{c}_2' : v > x \text{ and } \lambda(v) < \lambda(x) \text{ after 
$\C_1^+$-sorting}\}$ is nonempty.  However, the convexity of $\widetilde{P}$ 
makes elements of $S \setminus \bm{c}_1$ greater in $P$ than elements of 
$\bm{c}_2$, and in the same $\C_1^+$ chains as those elements.  Thus even if 
$\lambda(x)$ decreases, the label of $w$ is no larger.  Therefore the $\C_1^+$ 
chains remain sorted.  Similarly, the $\C_2^+$ chains remain sorted after 
first $\C_2^+$-sorting and then $\C_1^+$-sorting, so $P \in \n$ via 
$\{\C_1^+,\C_2^+\}$.  This coloring is induced by the coloring of 
$\widetilde{P} = \widetilde{P}' \cup \widetilde{w}$, since the coloring of 
$P'$ is induced by the coloring of $\widetilde{P}'$.  This concludes one 
direction of the proof.

Now let the $K$ element connected poset $P$ be in $\n$ via $\{\C_1,\C_2\}$.  It must be shown that $P$ reduces to a convex subposet of a cylinder poset, that it obeys the technical condition, and that $\C_1$ and $\C_2$ are induced by the coloring of this reduced poset.

Suppose there is a maximal or minimal element $z \in P$ incident to a doubly colored edge.  Let $P' = P \setminus z$.  Theorem~\ref{convexity} implies that $P' \in \n$ via $\{\C_1\vert_{P'},\C_2\vert_{P'}\}$.  By the induction hypothesis, $P' \leadsto \widetilde{P}'$ and the coloring of $P'$ is induced by the coloring of $\widetilde{P}'$.  In this situation, $P \leadsto \widetilde{P}'$ as well, and the technical condition for $P'$ together with Lemma~\ref{typeirules} for $P$ indicate that $P$ also satisfies the technical condition.  The chain covers $\C_i$ are induced by the coloring of $\widetilde{P}'$ because the chain covers $\C_i\vert_{P'}$ are as well, and the edge incident to $z$ is doubly colored.

It remains to consider when there is no such $z$.  As before, either $P$ is a chain, or there is a maximal or minimal element $w \in P$ adjacent to distinct elements, due to Lemma~\ref{chains}.  The case of a chain is straightforward, so suppose without loss of generality that a minimal element $w \in P$ is covered by $x$ and $y$.  Let $P' = P \setminus w$.  Theorem~\ref{convexity} indicates that $P' \in \n$ via $\{\C_1\vert_{P'},\C_2\vert_{P'}\}$, and each connected component of $P'$ has fewer than $K$ elements.  Thus, by the inductive hypothesis, $P' \leadsto \widetilde{P}'$, a convex subposet of a cylinder poset, $P'$ satisfies the technical condition, and the coloring of $P'$ is induced by that of $\widetilde{P}'$.  Let $P'$ be a convex subposet of $Cyl_{k,n}$.  Let $\widetilde{x}$ and $\widetilde{y}$ be the images of $x$ and $y$ in $\widetilde{P}'$.  If they are covered by a common element in $\widetilde{P}'$ (that is, $w$ is the minimum of a diamond in $P$), then $\widetilde{x}$ and $\widetilde{y}$ cover a common element $\widetilde{w} \in Cyl_{k,n}$.

Otherwise, removing $w$ either disconnects the poset or $\widetilde{P}'$ is 
connected and a convex subposet of the product of two chains.  Then $x$ and $y$ are each 
covered by at most one element in $P'$, and they are not covered by the same 
element.  If $\widetilde{x}$ is covered by $\widetilde{v} \in \widetilde{P}$, 
joined by a red edge, then the chain $[x,v]$ is red in $P'$ because of the 
induced coloring.  Thus the chain $[w,v]$ must be red in $P$ by 
Lemma~\ref{chains}.  If $x$ covers something other than $w$, then this edge 
must be blue in $\widetilde{P}'$.  Similar conclusions hold for the edges incident to $\widetilde{y}$, but the colors must be reversed since $w$ cannot be covered by two red edges in $P \in \n$.  Therefore the convexity of $\widetilde{P}'$ makes it possible either to draw the disjoint components of $\widetilde{P}'$ so that $\widetilde{x}$ and 
$\widetilde{y}$ both cover the same element $\widetilde{w}$ in $Cyl_{k,n}$, or, if $\widetilde{P}'$ has a single component, to 
choose $k$ and $n$ so that this is true.

As $\widetilde{w} \in Cyl_{k,n}$ is only covered by $\widetilde{x}$ and $\widetilde{y}$, the poset $\widetilde{P} = \widetilde{P}' \cup \widetilde{w}$ is convex in $Cyl_{k,n}$ and $P \leadsto \widetilde{P}$.  Lemma~\ref{typeirules} requires that the technical constraint be satisfied for $P$.  Finally, the coloring of $P$ is induced by that of $\widetilde{P}$ because $\{\C_1\vert_{P'},\C_2\vert_{P'}\}$ is induced by the coloring of $\widetilde{P}'$.
\end{proof}

The final case considered in the proof is when a maximal or minimal element of 
$P$ is adjacent to two other elements but is not in a diamond, and its removal 
does not disconnect the poset.  This describes a poset $P$ that 
can only reduce to a poset on the cylinder, while a maximal proper subposet of 
$P$ reduces to a convex subposet of the product of two chains.

\begin{figure}[htbp]
\begin{center}
$\begin{array}{c@{\hspace{.1in}}c}
\multicolumn{1}{l}{\mbox{\bf{(a)}}}
& \multicolumn{1}{l}{\mbox{\bf{(b)}}}\\
[-.45cm]
\epsfig{file=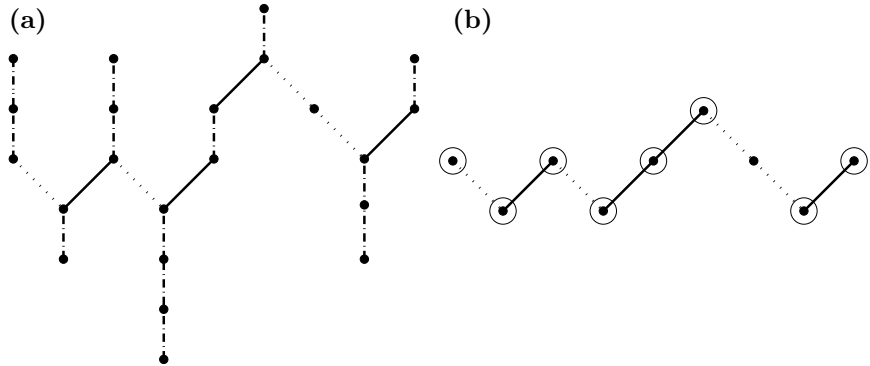} &
\raisebox{.73in}{\epsfig{file=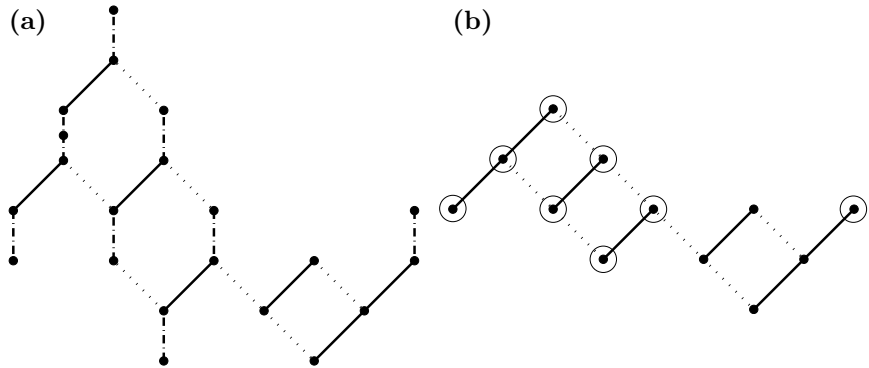}}
\end{array}$
\end{center}
\caption{$\bf{(a)}$ A poset $P \in \n$.
$\bf{(b)}$ The reduced poset $\widetilde{P}$.  Elements that split to form $P$ are circled.}
\label{treeanswer}
\end{figure}

\begin{figure}[htbp]
\begin{center}
$\begin{array}{c@{\hspace{.1in}}c}
\multicolumn{1}{l}{\mbox{\bf{(a)}}}
& \multicolumn{1}{l}{\mbox{\bf{(b)}}}\\
[-.45cm]
\epsfig{file=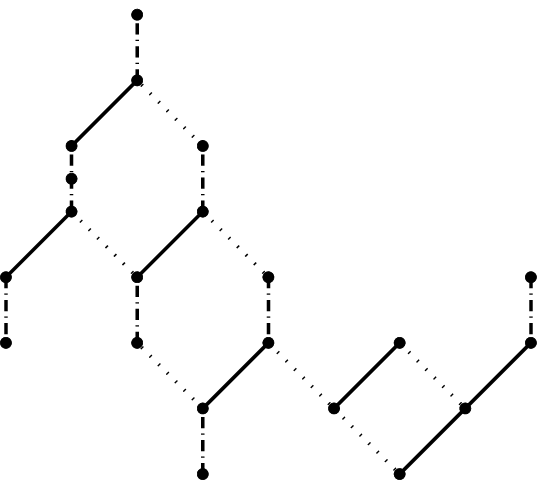} &
\raisebox{.27in}{\epsfig{file=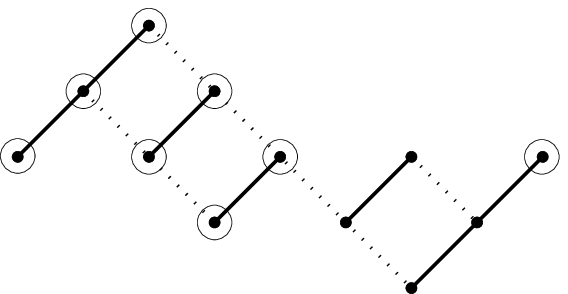}}
\end{array}$
\end{center}
\caption{$\bf{(a)}$ A poset $P \in \n$.
$\bf{(b)}$ The reduced poset $\widetilde{P}$.  Elements that split to form $P$ are circled.}
\label{bigcircuit}
\end{figure}

\begin{figure}[htbp]
\begin{center}
$\begin{array}{c@{\hspace{1in}}c}
\multicolumn{1}{l}{\mbox{\bf{(a)}}}
& \multicolumn{1}{l}{\mbox{\bf{(b)}}}\\
[-.45cm]
\hspace*{.35in}\epsfig{file=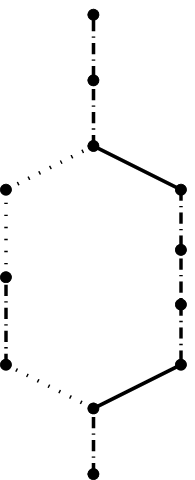} &
\hspace*{.35in}\raisebox{.55in}{\epsfig{file=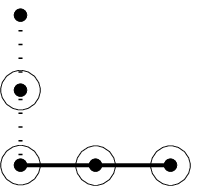}}
\end{array}$
\end{center}
\caption{$\bf{(a)}$ A poset $P \in \n$.
$\bf{(b)}$ A preimage of the reduced poset $\widetilde{P}$ in $Cyl_{2,5}$.  Elements that split to form $P$ are circled.}
\label{typeiistrings}
\end{figure}

\begin{figure}[htbp]
\begin{center}
$\begin{array}{c@{\hspace{.5in}}c}
\multicolumn{1}{l}{\mbox{\bf{(a)}}}
& \multicolumn{1}{l}{\mbox{\bf{(b)}}}\\
[-.45cm]
\hspace{.35in}\epsfig{file=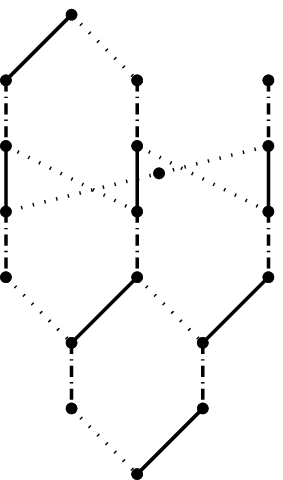} &
\hspace{.35in}\raisebox{.4in}{\epsfig{file=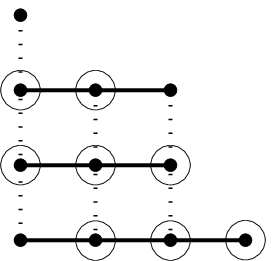}}
\end{array}$
\end{center}
\caption{$\bf{(a)}$ A poset $P \in \n$.
$\bf{(b)}$ A preimage of the reduced poset $\widetilde{P}$ in $Cyl_{3,7}$.  Elements that split to form $P$ are circled.}
\label{bigcrown}
\end{figure}

Examples of non-messing-up posets are depicted in Figures~\ref{treeanswer}(a), \ref{bigcircuit}(a), \ref{typeiistrings}(a), and~\ref{bigcrown}(a).  The first two of these reduce to convex subposets of the product of two chains, while the last two do not.

Because of the definition of an induced coloring, a Type II diamond as described in Lemma~\ref{diamondcolors} occurs only when a diamond a non-messing-up poset goes around the cylinder.  This explains the technical condition.

The requirement for membership in $\n$ is the existence of a pair of chain
covers $\{\C_1, \C_2\}$ with particular properties.  It is natural to ask if there are other choices for the chain covers.  A poset of the form depicted in Figure~\ref{typeiistrings}(a), that is, a poset consisting of a single diamond and its bottom and top chains, can also be colored so that the diamond has Type I if the bounds of Lemma~\ref{typeirules}
are satisfied.  Otherwise, the only freedom in defining the chain covers arises
from the various ways to reduce $P$ due to splits as depicted in
Figure~\ref{splits}(a).

\section{Further directions}

The classification of $\n$ prompts further questions relating to the non-messing-up property.  The final section of this paper suggests several such questions and provide answers to some.

\subsection{The set $\n ' \subsetneq \n$ with reduced redundancy}\label{n'}\

In the classification of $\n$, there were instances of a $\C_i$
chain entirely contained in a $\C_{3-i}$ chain.  These chain covers have the non-messing-up property, but there is a certain redundancy: this particular $\C_i$ chain
adds no information about the relations in the poset since
its labels are already ordered after the $\C_{3-i}$-sort.

\begin{defn}
The class $\n'$ consists of all posets $P \in \n$ via $\{\C_1,\C_2\}$ such that $\bm{c}_i \not \subseteq \bm{c}_{3-i}$
for all $\bm{c}_1 \in \C_1$ and $\bm{c}_2 \in \C_2$.
\end{defn}

Because the coloring of a non-messing-up poset is induced by its reduced poset, the elements of $\n'$ can be determined by looking at these reduced posets.  Consider some $P \in \n$ that reduces to $\widetilde{P}$, a convex
subposet of a cylinder poset, where $\widetilde{P} \in \n$ via $\{\widetilde{\C}_1,\widetilde{\C}_2\}$.  Suppose there is an element $\widetilde{v} \in
\widetilde{P}$ that is adjacent to at most one other element in
$\widetilde{P}$.  Then $\widetilde{v} \in \widetilde{\bm{c}}_i \in
\widetilde{\C}_i$ such that, up to color reversal, $\widetilde{\bm{c}}_1 = \{v\}$, and consequently $\widetilde{\bm{c}}_1
\subseteq \widetilde{\bm{c}}_2$.  However elements of
$\widetilde{P}$ split to form $P$, the resulting chain $\bm{c}_1$ will be
entirely contained in the resulting chain $\bm{c}_2$.  Therefore,
$P \notin \n'$.

Consider the other elements in $\widetilde{P}$.  No red chain that intersects a Type I diamond in an edge will be contained in a blue chain, and likewise with the colors reversed.  Consider a chain in $\widetilde{P}$ that shares no edge with any Type I diamonds.  Call this a \emph{branch chain}, and a maximal such chain a \emph{maximal branch chain}.  Let $\widetilde{P}$ be the most reduced version of $P$, so every maximal branch chain consists of two or three elements in $\widetilde{P}$.  The chain covers for $\widetilde{P}$ induce the chain covers for $P$, which yields the following result.

\begin{thm}\label{n'answer}
The collection $\n'$ is the set of posets in $\n$ where every maximal branch chain in the reduced poset $\widetilde{P}$ consists of exactly two elements, and every element of $\widetilde{P}$ is adjacent to at least two other elements in $\widetilde{P}$.
\end{thm}

If an element of $\n'$ consists of a single diamond and its top and bottom chains, then that diamond must have Type I.  Similarly, there are no posets whose Hasse diagrams are trees, as in Figure~\ref{treeanswer}(a), as every reduced tree has an element with at most one incident edge.

\subsection{The set $\n'' \subseteq \n$ with reduced redundancy}\

This section considers another notion of redundancy for elements of $\n$.  In Theorem~\ref{nmu}, the rows and
columns have minimal redundancy in the sense that for any
row $\bm{r}$ and any column $\bm{c}$, $\# (\bm{r} \cap \bm{c}) = 1$.  

\begin{defn}
The class $\n''$ consists of all posets $P \in \n$ via $\{\C_1,\C_2\}$ such that $\#(\bm{c}_1 \cap \bm{c}_2) \le 1$ for all $\bm{c}_i \in \C_i$.
\end{defn}

Observe that $\n' \not\subseteq \n''$ and $\n'' \not\subseteq \n'$, since
$\n''$ permits a single element chain in $\C_i$,
necessarily contained in a chain of $\C_{3-i}$, and elements of $\n'$ can have chain intersections of any size.

The classification of $\n$ is based on two classes of allowable posets and
the posets that result from splitting elements of these in particular
ways.  Consider the least reduced poset for an element of $\n$.  This is a convex subposet of a cylinder poset, where the splits are of the smallest size.  When elements are split, the new edges must be doubly colored by Lemma~\ref{chains}.  Thus the
only connected posets that can be in $\n''$ are themselves convex subposets of
a cylinder poset.  The coloring inherited from the posets described in
Theorem~\ref{nmu} and Theorem~\ref{cylinder} colors each of
these posets so that any two differently colored chains intersect
in at most one element.

\begin{thm}
The collection $\n''$ is the set of finite posets each of whose connected components is a convex subposet of a cylinder poset.
\end{thm}

\subsection{Open questions}\

This paper generalizes Theorem~\ref{nmu} by characterizing $\n$, and characterizes the more restrictive classes $\n'$ and $\n''$.  It may also be fruitful to examine other generalizations and related topics, some of which are suggested here.

This paper studies finite posets and saturated chains, but interesting
results may arise if one or both of these restrictions are relaxed.  
Similarly, one could study posets with some variation of the non-messing-up
phenomenon.  For example, one could consider more than two sets of chains,
or expand beyond identities like $S_iS_{3-i}S_i(\mathcal{L}(P)) =
S_{3-i}S_i(\mathcal{L}(P))$ for all labelings $\mathcal{L}$ of $P$ and $i
\in \{1,2\}$, where $S_i(\mathcal{L}(P))$ represents $\C_i$-sorting a
labeling $\mathcal{L}$ of a poset $P$.

Additionally, as stated earlier, any labeling of a poset $P \in \n$ produces a linear extension of $P$ after performing the two sorts.  It would be interesting to understand the distribution of the linear extensions that arise in this way.

These are examples of issues related to the non-messing-up phenomenon that warrant further study.  The author hopes to address some of them in the future.

\end{document}